\documentclass[12pt]{article} 
\usepackage{amssymb,amsmath,color}
\usepackage{amsthm}
\input{psfig.sty}

\newtheorem{theorem}{Theorem}[section]
\newtheorem{proposition}[theorem]{Proposition}
\newtheorem{lemma}[theorem]{Lemma}
\newtheorem{corollary}[theorem]{Corollary}
\newtheorem{conjecture}[theorem]{Conjecture}

\theoremstyle{definition}

\theoremstyle{remark}

\begin{document}
\newcommand{\beq}{\begin{equation}} 
\newcommand{\eeq}{\end{equation}}
\newcommand{\zz}{\mathbb{Z}}
\newcommand{\pp}{\mathbb{P}} 
\newcommand{\qq}{\mathbb{Q}} 
\newcommand{\nn}{\mathbb{N}}
\newcommand{\rr}{\mathbb{R}}
\newcommand{\bm}[1]{{\mbox{\boldmath $#1$}}}
\newcommand{\con}{\mathrm{Comp}(n)}
\newcommand{\sn}{\mathfrak{S}_n} 
\newcommand{\fs}{\mathfrak{S}}
\newcommand{\st}{\,:\,} 
\newcommand{\ep}{$\ \ \Box$}
\newcommand{\as}{\mathrm{as}}
\newcommand{\is}{\mathrm{is}}
\newcommand{\lgn}{\mathrm{len}}
\newcommand{\op}{\mathrm{OP}}
\newcommand{\rsk}{\stackrel{\mathrm{rsk}}{\rightarrow}}
\newcommand{\dis}{\displaystyle}
\newcommand{\co}{\mathrm{co}}
\newcommand{\cp}{\mathrm{Comp}}
\newcommand{\alt}{\mathrm{alt}}
\newcommand{\ralt}{\mathrm{ralt}}
\newcommand{\bea}{\begin{eqnarray}}
\newcommand{\twoline}[2]{\genfrac{}{}{0pt}{}{#1}{#2}}
\newcommand{\at}{\tan^{-1}}

\newcommand{\beas}{\begin{eqnarray*}}
\newcommand{\eea}{\end{eqnarray}}
\newcommand{\eeas}{\end{eqnarray*}}
\newcommand{\odd}{\,\mathrm{odd}}
\newcommand{\even}{\,\mathrm{even}}
\newcommand{\lm}{\lambda/\mu}

\thispagestyle{empty}

\vskip 20pt
\begin{center}
\textcolor{red}{{\large\bf Alternating Permutations and Symmetric
    Functions}} 
\vskip 15pt
{\textcolor{blue}{\bf Richard P. Stanley}}\\
{\it Department of Mathematics, Massachusetts Institute of
Technology}\\
{\it Cambridge, MA 02139, USA}\\
{\texttt{rstan@math.mit.edu}}\\[.2in]
{\bf\small version of 21 June 2006}\\
\end{center}

\begin{abstract}
We use the theory of symmetric functions to enumerate various classes
of alternating permutations $w$ of $\{1,2,\dots,n\}$. These classes
include the following: (1) both $w$ and $w^{-1}$ are alternating, (2)
$w$ has certain special shapes, such as $(m-1,m-2,\dots,1)$, under the
RSK algorithm, (3) $w$ has a specified cycle type, and (4) $w$ has a
specified number of fixed points. We also
enumerate alternating permutations of a multiset. Most of our
formulas are umbral expressions where after expanding the expression
in powers of a variable $E$, $E^k$ is interpreted as the Euler
number $E_k$. As a small corollary, we obtain a combinatorial
interpretation of the coefficients of an asymptotic expansion
appearing in Ramanujan's Lost Notebook.
\end{abstract}

\section{Introduction.} \label{sec1}
\indent This paper can be regarded as a sequel to the classic paper
\cite{foulkes} of H. O. Foulkes in which he relates the enumeration of
alternating permutations to the representation theory of the symmetric
group and the theory of symmetric functions. We assume familiarity
with symmetric functions as presented in \cite[Ch.~7]{ec2}.
Let $\sn$ denote the symmetric group of all permutations of
$1,2,\dots,n$. A permutation $w=a_1 a_2\cdots a_n\in\sn$ is
\emph{alternating} if $a_1>a_2<a_3>a_4<\cdots$. Equivalently, write
$[m]=\{1,2,\dots,m\}$ and define the \emph{descent set} $D(w)$ of
$w\in\sn$ by
  $$ D(w)  = \{ i\in[n-1]\st a_i>a_{i+1}\}. $$
Then $w$ is alternating if $D(w)=\{1,3,5,\dots\}\cap [n-1]$.
Similarly, define $w$ to be \emph{reverse alternating} if
$a_1<a_2>a_3<a_4>\cdots$.  Thus $w$ is reverse alternating if
$D(w)=\{2,4,6,\dots\}\cap [n-1]$. Also define the \emph{descent
  composition} $\co(w)$ by 
  \beq \co(w) = (\alpha_1,\alpha_2,\dots,\alpha_k), 
    \label{eq:cow} \eeq
where $D(w)=\{\alpha_1,\alpha_1+\alpha_2,\dots,\alpha_1+\cdots+
\alpha_{k-1}\}$ and $\sum \alpha_i=n$. Thus $\alpha\in\cp(n)$, where
$\cp(n)$ denotes the set of compositions of $n$. 

Let $E_n$ denote the number of alternating permutations in $\sn$. Then
$E_n$ is called an \emph{Euler number} and was shown by D. Andr\'e
\cite{andre} to satisfy
  \beq \sum_{n\geq 0} E_n\frac{x^n}{n!} = \sec x + \tan x. 
   \label{eq:eudef} \eeq
(Sometimes one defines $\sum (-1)^nE_nx^{2n}/(2n)! = \sec x$, but we
will adhere to (\ref{eq:eudef}).)  Thus $E_{2m}$ is also called a
\emph{secant number} and $E_{2m+1}$ a \emph{tangent number}. The
bijection $w\mapsto w'$ on $\sn$ defined by $w'(i)=n+1-w(i)$ shows
that $E_n$ is also the number of reverse alternating permutations in
$\sn$. However, for some of the classes of permutations considered
below, alternating and reverse alternating permutations are not
equinumerous.

Foulkes defines a certain (reducible) representation of $\sn$ whose
dimension is $E_n$. He shows how this result can be used to compute
$E_n$ and other numbers related to alternating permutations, notably
the number of $w\in\sn$ such that both $w$ and $w^{-1}$ are
alternating. Foulkes' formulas do not give a ``useful'' computational
method since they involve sums over partitions whose terms involve
Littlewood-Richardson coefficients. We show how Foulkes' results can
actually be converted into useful generating functions for computing
such numbers as (a) the number of alternating permutations $w\in\sn$
with conditions on their cycle type (or conjugacy class). The special
case of enumerating alternating involutions was first raised by
Ehrenborg and Readdy and discussed further by Zeilberger
\cite{zeil}. Another special case is that of alternating permutations
with a specified number of fixed points. Our proofs use, in addition
to Foulkes' representation, a result of Gessel and Reutenauer
\cite{g-r} on permutations with given descent set and cycle type. (b)
The number of $w\in\sn$ such that both $w$ and $w^{-1}$ are
alternating, or such that $w$ is alternating and $w^{-1}$ is reverse
alternating. (c) The number of alternating permutations of certain
shapes (under the RSK algorithm). (d) The number of alternating
permutations of a multiset of integers, under various interpretations
of the term ``alternating.''

\textsc{Acknowledgment.} I am grateful to Ira Gessel for providing
some useful background information and references and to an anonymous
referee for several helpful comments, in particular, pointing out a
gap in the proof of Corollary~\ref{cor:asy}.

\section{The work of Foulkes.} \label{sec2}
We now review the results of Foulkes that will be the basis for our
work. Given a composition $\alpha$ of $n$, let $B_\alpha$ denote the
corresponding border strip (or ribbon or skew hook) shape as defined
e.g.\ in \cite{foulkes2}\cite[p.~383]{ec2}. Let $s_{B_\alpha}$ denote
the skew Schur function of shape $B_\alpha$. The following result of
Foulkes \cite[Thm.~6.2]{foulkes2} also appears in
\cite[Cor.~7.23.8]{ec2}.

\begin{theorem} \label{thm:f1}
Let $\alpha$ and $\beta$ be compositions of $n$. Then
  $$ \langle s_{B_\alpha},s_{B_\beta}\rangle = \#\{w\in\sn\st
     \co(w)=\beta,\ \co(w^{-1})=\alpha\}. $$
\end{theorem}

We let $\tau_n=B_\alpha$ where $\alpha=(1,2,2,\dots,2,j)\in\cp(n)$,
where $j=1$ if $n$ is even and $j=2$ if $n$ is odd. Thus if $'$
indicates conjugation (reflection of the shape about the main
diagonal), then $\tau'_{2k+1}=\tau_{2k+1}$, while $\tau'_{2k} =
(2,2,\dots, 2)$. We want to expand the skew Schur functions
$s_{\tau_n}$ and $s_{\tau'_n}$ in terms of power sum symmetric
functions. For any skew shape $\lambda/\mu$ with $n$ squares, let
$\chi^{\lambda/\mu}$ denote the character of $\sn$ satisfying
ch$(\chi^{\lambda/\mu})=s_{\lambda/\mu}$. Thus by the definition
\cite[p.~351]{ec2} of ch we have
  $$ s_{\lambda/\mu} = \sum_{\rho\vdash n} z_\rho^{-1}
  \chi^{\lambda/\mu}(\rho)p_\rho, $$ 
where $\chi^{\lambda/\mu}(\rho)$ denotes the value of
$\chi^{\lambda/\rho}$ at any permutation $w\in\sn$ of cycle type
$\rho$. 

The main result \cite[Thm.~6.1]{foulkes}\cite[Exer.~7.64]{ec2}
of Foulkes on the connection between alternating permutations and
representation theory is the following.

\begin{theorem} \label{thm:foulkes}
 \emph{(a)} Let $\mu\vdash n$, where $n=2k+1$. Then
  $$ \chi^{\tau_n}(\mu) = \chi^{\tau'_n}(\mu) = \left\{ 
    \begin{array}{rl} 0, & \mbox{if $\mu$ has an even part}\\[.05in]
         (-1)^{k+r}E_{2r+1}, & \mbox{if $\mu$ has $2r+1$ odd parts
         and}\\  & \mbox{\ \ no even parts}. \end{array} \right. $$
 \emph{(b)} Let $\mu\vdash n$, where $n=2k$. Suppose that $\mu$ has
$2r$ odd parts and $e$ even parts. Then
   \beas \chi^{\tau_n}(\mu) & = &  (-1)^{k+r+e}E_{2r}\\[.05in]
     \chi^{\tau'_n}(\mu) & = &  (-1)^{k+r}E_{2r}. \eeas
\end{theorem}

\textsc{Note.} Foulkes obtains his result from the Murnaghan-Nakayama
rule. It can also be also be obtained from the formula 
  $$ \sum_{n\geq 0} s_{\tau_n}t^n = \frac{1}{\sum_{n\geq 0}(-1)^n
      h_{2n}t^{2n}}+\frac{\sum_{n\geq 0}h_{2n+1}t^{2n+1}}
      {\sum_{n\geq 0}(-1)^mh_{2n}t^{2n}}, $$
where $s_{\tau_n}$ denotes a skew Schur function. This formula is due
to  Carlitz \cite{carlitz} and is also stated at the bottom of
page~520 of \cite{ec2}. 

Foulkes' result leads immediately to our main tool in what follows.
Throughout this paper we will use \emph{umbral notation} \cite{r-t}
for Euler numbers. In other words, any polynomial in $E$ is to be
expanded in terms of powers of $E$, and then $E^k$ is replaced by
$E_k$. The replacement of $E^k$ by $E_k$ is always the \emph{last} step in
the evaluation of an umbral expression.  For instance,
  $$ (E^2-1)^2 = E^4 -2E^2+1=E_4-2E_2+1=5-2\cdot 1+1=4. $$
Similarly,
  \beas (1+t)^E & = & 1+Et+{E\choose 2}t^2+{E\choose 3}t^3+\cdots\\
    & = & 1+Et+\frac 12(E^2-E)t^2+\frac 16(E^3-3E^2+2E)t^3+\cdots\\
    & = & 1+Et+\frac 12(E_2-E_1)t^2+\frac 16(E_3-3E_2+2E_1)t^3+
       \cdots\\ & = & 
     1+1\cdot t+\frac 12(1-1)t^2+\frac 16(2-3\cdot 1+2\cdot 1)t^3
        + \cdots\\ & = & 1+t+\frac 16 t^3+\cdots. \eeas
If $f=f(x_1,x_2,\dots)$ is a
symmetric function then we use the notation $f[p_1,p_2,\dots]$ for $f$
regarded as a polynomial in the power sums. For instance, if
$f=e_2=\sum_{i<j}x_i x_j=\frac 12(p_1^2-p_2)$ then
  $$ e_2[E,-E,\dots]= \frac 12(E^2+E)=1. $$ 

\begin{theorem} \label{thm:main}
Let $f$ be a homogenous symmetric function of degree $n$. If $n$ is
  odd then 
  \beq \langle f,s_{\tau_n}\rangle =\langle f,s_{\tau'_n}\rangle =
         f[E,0,-E,0,E,0,-E,\dots] \label{eq:feo} \eeq
If $n$ is even then 
 \beas \langle f,s_{\tau_n}\rangle & = & 
         f[E,-1,-E,1,E,-1,-E,1,\dots]\\
   \langle f,s_{\tau'_n}\rangle & = & 
         f[E,1,-E,-1,E,1,-E,-1,\dots]. \eeas
\end{theorem}

\textbf{Proof.} Suppose that $n=2k+1$. Let $\op(n)$ denote the set
of all partitions of $n$ into odd parts. If $\mu\in\op(n)$ and $\mu$
has $\ell(\mu)=2r+1$ (odd) parts, then write $r=r(\mu)$. Let
$f=\sum_{\lambda\vdash n} c_\lambda p_\lambda$. Then by
Theorem~\ref{thm:foulkes} we have 
  \beas \langle f,s_{\tau_n}\rangle & = & \left\langle
   \sum_\lambda c_\lambda p_\lambda, \sum_{\mu\in\op(n)}z_\mu^{-1}
    (-1)^{k+r(\mu)} E_{\ell(\mu)}p_\mu\right\rangle\\ & = &
   \sum_{\mu\in\op(n)} c_\mu (-1)^{k+r(\mu)}E_{\ell(\mu)}. \eeas
If $\mu\in \op(n)$ and we substitute $(-1)^jE$ for $p_{2j+1}$ in
$p_\mu$ then we obtain 
  \beas \prod_{i=1}^{\ell(\mu)}(-1)^{\frac 12(\mu_i-1)}E & = &
          (-1)^{\frac 12(2k+1-(2r(\mu)+1))}E^{\ell(\mu)}\\
   & = & (-1)^{k+r(\mu)}E^{\ell(\mu)}, \eeas
and equation~(\ref{eq:feo}) follows. The case of $n$ even is
analogous. \qed 

\section{Inverses of alternating permutations.}
In this section we derive generating functions for the number of
alternating permutations in $\sn$ whose inverses are alternating or
reverse alternating. This problem was considered by Foulkes
\cite[{\S}5]{foulkes}, but his answer does not lend itself to easy
computation. Such ``doubly alternating'' permutations were also
considered by Ouchterlony \cite{ouch} in the setting of pattern
avoidance. A special class of doubly alternating permutations, viz.,
those that are Baxter permutations, were enumerated by Guibert and
Linusson \cite{g-l}.

\begin{theorem} \label{thm:doubalt}
Let $f(n)$ denote the number of permutations $w\in\sn$ such that both
$w$ and $w^{-1}$ are alternating, and let $f^\ast(n)$ denote the
number of $w\in\sn$ such that $w$ is alternating and $w^{-1}$ is
reverse alternating. Let 
    \beas L(t) & = & \frac 12\log\frac{1+t}{1-t}\\ & = &
    t+\frac{t^3}{3}+\frac{t^5}{5}+\cdots. \eeas
Then
    \bea \sum_{k\geq 0}f(2k+1)t^{2k+1} & = & \sum_{r\geq 0}E_{2r+1}^2
    \frac{L(t)^{2r+1}}{(2r+1)!} \label{eq:aa1}\\[.1in]
    f^\ast(2k+1) & = & f(2k+1) \label{eq:aa2}\\[.1in]
   \sum_{k\geq 0}f(2k)t^{2k} & = & \frac{1}{\sqrt{1-t^2}} \sum_{r\geq
     0}E_{2r}^2 \frac{L(t)^{2r}}{(2r)!} \label{eq:aa3}\\[.1in]
    f^\ast(2k) & = & f(2k)-f(2k-2). \label{eq:aa4} \eea
\end{theorem}

\textbf{Proof.} By Theorem~\ref{thm:f1} we have $f(n) =\langle
s_{\tau_n}, s_{\tau_n}\rangle$. Let $n=2k+1$. Then it follows from
Theorems~\ref{thm:foulkes} and \ref{thm:main} that (writing
$r=r(\mu)$) 
  \bea f(n) & = & \sum_{\mu\in\op(n)}
  z_\mu^{-1}(-1)^{k+r}E_{2r+1}(-1)^{k+r} E^{2r+1}
  \nonumber \\ & = &
    \sum_{\mu\in\op(n)} z_\mu^{-1}E_{2r+1}^2. \label{eq:fnodd} \eea
Now by standard properties of exponential generating functions
\cite[{\S}5.1]{ec2} or by specializing the basic identity 
  $$ \sum_\lambda z_\lambda^{-1}p_\lambda=\exp \sum_{n\geq 1}
     \frac 1np_n, $$
we have
  \beas \sum_{k\geq 0}\sum_{\mu\in\op(2k+1)}
  z_\mu^{-1}y^{\ell(\mu)}t^{2k+1} & = & 
   \exp\left(y\left(t+\frac{t^3}{3}+\frac{t^5}{5}+\cdots\right)\right)
   \\ & = & \exp( yL(t)). \eeas
The coefficient of $y^{2r+1}$ in the above generating function is
therefore $L(t)^{2r+1}/(2r+1)!$, and the proof of (\ref{eq:aa1})
follows.

Since $\tau_n=\tau'_n$ for $n$ odd we have
  $$ f^\ast(n) =\langle s_{\tau_n},s_{\tau'_n}\rangle =
   \langle s_{\tau_n},s_{\tau_n}\rangle = f(n), $$
so (\ref{eq:aa2}) follows.

The argument for $n=2k$ is similar. For $\mu\vdash n$ let $e=e(\mu)$
  denote the number of even parts of $\mu$ and $2r=2r(\mu)$ the number
  of odd parts. Now the relevant formulas for computing $f(n)$ are
  \beas f(n) & = & \sum_{\mu\vdash n}
  z_\mu^{-1}(-1)^{k+r+e}E_{2r}(-1)^{k+r+e} E^{2r}\\ & = &
    \sum_{\mu\vdash n} z_\mu^{-1}E_{2r}^2 \eeas
and
  \beas \sum_{k\geq 0}\sum_{\mu\vdash n}
  z_\mu^{-1}y^{2r(\mu)}t^n & = & 
   \exp\left(y\left(t+\frac{t^3}{3}+\frac{t^5}{5}+\cdots\right)+
    \left( \frac{t^2}{2}+\frac{t^4}{4}+\cdots\right)\right)\\
   & = & (1-t^2)^{-1/2}(\exp( yL(t)), \eeas
from which (\ref{eq:aa3}) follows.

For the case $f^\ast(n)$ when $n$ is even we have
  \beas f^\ast(n) & = & \sum_{\mu\vdash n}
  z_\mu^{-1}(-1)^{k+r+e}E_{2r}(-1)^{k+r} E^{2r}\\ & = &
    \sum_{\mu\vdash n} z_\mu^{-1}(-1)^eE_{2r}^2 \eeas
and
  \beas \sum_{k\geq 0}\sum_{\mu\vdash 2k}
  z_\mu^{-1}(-1)^{e(\mu)}y^{2r(\mu)}t^{2k} & = & 
   \exp\left(y\left(t+\frac{t^3}{3}+\frac{t^5}{5}+\cdots
   \right)\right.\\  & & \quad 
    -\left.\left( \frac{t^2}{2}+\frac{t^4}{4}+
  \cdots\right)\right)\\[.1in]
   & = & \sqrt{1-t^2}\exp yL(t). \eeas
Hence
   $$ \sum_{k\geq 0}f^\ast(2k)t^{2k} = 
     (1-t^2)\sum_{k\geq 0}f(2k)t^{2k}, $$
from which (\ref{eq:aa4}) follows. \qed

Whenever we have explicit formulas or generating functions for
combinatorial objects we can ask for combinatorial proofs of
them. Bruce Sagan has pointed out that equation (\ref{eq:aa2}) follows
from reversing the permutation, i.e., changing $a_1 a_2 \cdots a_n$ to
$a_n \cdots a_2 a_1$.  We do not know combinatorial proofs of
equations (\ref{eq:aa1}), (\ref{eq:aa3}) and (\ref{eq:aa4}). To prove
equations (\ref{eq:aa1}) and (\ref{eq:aa3}) combinatorially, we
probably need to interpret them as exponential generating functions,
e.g., write the left-hand side of (\ref{eq:aa1}) as $\sum_{k\geq 0}
(2k+1)!f(2k+1)t^{2k+1}/(2k+1)!$. Let us also note that if $g(n)$
denotes the number of reverse alternating $w\in\sn$ such that $w^{-1}$
is also reverse alternating, then $f(n)=g(n)$ for all $n$. This fact
can be easily shown using the proof method above, and it is also a
consequence of the RSK algorithm. For suppose that $w$ and $w^{-1}$
are alternating, $w\rsk (P,Q)$ and $w'\rsk (P^t,Q^t)$ (where $^t$
denotes transpose). Then by \cite[Lemma~7.23.1]{ec2} the map $w\mapsto
w'$ is a bijection between permutations $w\in\sn$ such that both $w$
and $w^{-1}$ are alternating, and permutations $w'\in\sn$ such that
both $w$ and $(w')^{-1}$ are reverse alternating.  Is there a simpler
proof that $f(n)=g(n)$ avoiding RSK?

\section{Alternating tableaux of fixed shape.}
Let $T$ be a standard Young tableau (SYT). The \emph{descent set}
$D(T)$ is defined by \cite[p.~351]{ec2}
  $$ D(T) = \{i\st i+1\ \mbox{is in a lower row than}\ i\}. $$
For instance, if
    $$ T = \begin{array}{l}1\,2\,5\\ 3\,4\\ 6, \end{array} $$ 
then $D(T)=\{2,5\}$. We also define the \emph{descent composition}
$\co(T)$ in analogy with equation~(\ref{eq:cow}).  A basic property of
the RSK algorithm asserts that $D(w)=D(Q)$ if $w\rsk (P,Q)$.
An SYT $T$ of size $n$ is called \emph{alternating} if
$D(T)=\{1,3,5,\dots\}\cap[n-1]$ and \emph{reverse alternating} if
$D(T)=\{2,4,6,\dots\}\cap[n-1]$. The following result is an immediate
consequence of Theorem~7.19.7 and Corollary~7.23.6 of \cite{ec2}.
 
\begin{theorem} \label{thm:dco} 
Let $\lambda\vdash n$ and $\alpha\in\con$. Then $\langle
s_\lambda,s_{B_\alpha}\rangle$ is equal to the number of SYT of shape
$\lambda$ and descent composition $\alpha$.
\end{theorem}

Let $\alt(\lambda)$ (respectively, $\ralt(\lambda)$) denote the number
of alternating (respectively, reverse alternating) SYT of shape
$\lambda$. The following result then follows from
Theorems~\ref{thm:main} and \ref{thm:dco}.

\begin{theorem} \label{thm:shape}
Let $\lambda\vdash n$ and $\alpha\in\con$. If $n$ is odd, then 
  $$ \alt(\lambda) = \ralt(\lambda)
      =s_\lambda[E,0,-E,0,E,0,-E,\dots]. $$ 
If $n$ is even then 
   \beas \alt(\lambda) & = & s_\lambda[E,-1,-E,1,E,-1,-E,1,\dots]\\ 
     \ralt(\lambda) & = & s_\lambda[E,1,-E,-1,E,1,-E,-1,\dots].
   \eeas
\end{theorem}
 
Theorem~\ref{thm:shape} ``determines'' the number of alternating SYT
of any shape $\lambda$, but the formula is not very enlightening. We
can ask whether there are special cases for which the formula can be
made more explicit. The simplest such case occurs when $\lambda$ is
the ``staircase'' $\delta_m=(m-1,m-2,\dots,1)$. For any partition
$\lambda$ write $H_\lambda$ for the product of the hook lengths of
$\lambda$ \cite[p.~373]{ec2}. For instance,
  $$ H_{\delta_m}= 1^{m-1}\,3^{m-2}\,5^{m-3} \cdots (2m-3). $$

\begin{theorem} \label{thm:stair}
If $m=2k$ then
  $$ \alt(\delta_m)=\ralt(\delta_m) =
       E^k\prod_{j=1}^{m-2} (E^2+j^2)^{k-\lceil j/2 \rceil}. $$ 
If $m=2k+1$ then
  $$ \alt(\delta_m)=\ralt(\delta_m) =
   E^k\prod_{j=1}^{m-2} (E^2+j^2)^{k-\lfloor j/2 \rfloor}. $$ 
\end{theorem}

\textbf{Proof.} By the Murnaghan-Nakayama rule, $s_{\delta_m}$ is a
polynomial in the odd power sums $p_1,p_3,\dots$
\cite[Prop.~7.17.7]{ec2}. Assume that $m$ is odd. Then by the
hook-content formula \cite[Cor.~7.21.4]{ec2} we have 
  \bea s_{\delta_m}[E,0,E,0,\dots] & = & s_{\delta_m}[E,E,E,\dots]
   \nonumber \\ & = & 
   \frac{E^k\prod_{j=1}^{m-2} (E^2-j^2)^{k-\lfloor j/2 \rfloor}}
   {H_{\delta_m}}.
     \label{eq:seee} \eea
Let $n={m\choose 2}$, and suppose that $n$ is odd, say $n=2r+1$.
Let $\lambda\in\op_n$ and $2j+1=\ell(\lambda)$. Thus
  $$ p_\lambda[E,0,E,0,E,0,\dots]=E^{2j+1}. $$
A simple parity argument shows that  
    $$ p_\lambda[E,0,-E,0,E,0,-E,0,\dots] = (-1)^{r-j}E^{2j+1}. $$ 
It follows that we obtain $s_{\delta_m}[E,0,-E,0,E,0,-E,0,\dots]$ from
the polynomial expansion of $s_{\delta_m}[E,0,E,0,E,0,\dots]$ by
replacing each power $E^{2j+1}$ with $(-1)^{r-j}E^{2j+1}$. The proof
for $m$ odd and $n$ odd now follows from equation~(\ref{eq:seee}). 

The argument for the remaining cases, viz., (a) $m$ odd, $n$ even, (b)
$m$ even, $n$ odd, and (c) $m$ even, $n$ even, is completely
analogous. \qed

There are some additional partitions $\lambda$ for which
$\alt(\lambda)$ and $\ralt(\lambda)$ factor nicely as polynomials in
$E$. One such case is the following.

\begin{theorem} \label{thm:opsdet}
Let $p$ be odd, and let $p\times p$ denote the partition of $p^2$
whose shape is a $p\times p$ square. Then
  \beas \alt(p\times p) & = & \ralt(p\times p)\\ & = &
    \frac{E^p(E^2+2^2)^{p-1}(E^2+4^2)^{p-2}\cdots (E^2+(2(p-1))^2)}
       {H_{p\times p}}. \eeas
\end{theorem}

\emph{Proof} (sketch). Let $h_n$ denote the complete symmetric
function of degree $n$. From the identity
  $$ \sum_{n\geq 0} h_nt^n =\exp\sum_{n\geq 1}\frac{p_nt^n}{n} $$
we obtain
  \beas \sum_{n\geq 0}h_n[E,0,E,0,E,0,\dots]t^n & = & 
     \exp \sum_{n\odd}\frac{Et^n}{n}\\ & = &
       \left(\frac{1+t}{1-t}\right)^{E/2}. \eeas 
Write
  $$ \left(\frac{1+t}{1-t}\right)^{E/2}= \sum_{n\geq 0} a_n(E)t^n. $$
The Jacobi-Trudi identity \cite[{\S}7.16]{ec2} implies that
$s_{p\times p} = \det(h_{p-i+j})_{i,j=1}^p$. Hence
  \beq s_{p\times p}[E,0,E,0,E,0,\dots] =
    \det(a_{p-i+j}(E))_{i,j=1}^p. \label{eq:oddps} \eeq
I am grateful to Christian Krattenthaler and Dennis Stanton for
evaluating the above determinant. Krattenthaler's argument is as
follows. Write 
  $$ \left(\frac{1+t}{1-t}\right)^{E/2} =\left( 1+\frac{2t}{1-t}
    \right)^{E/2} =1+\sum_{n\geq 1}t^n\sum_{k=1}^n{n-1\choose k-1}
    {E/2\choose k}2^k. $$
After substituting $k+1$ for $k$, we see that we want to compute the
Hankel determinant 
  $$ \det_{0\leq i,j\leq n} \left( \sum_{k=0}^{i+j}
    {i+j\choose k}{E/2\choose k+1}2^{k+1}\right). $$
Now by a folklore result \cite[Lemma~15]{kratt} we conclude that this
determinant is the same as 
  $$ \det_{0\leq i,j\leq n} \left( {E/2\choose i+j+1}2^{i+j+1}
        \right). $$
When this determinant is expanded all powers of 2 are the same, so we
are left with evaluating
   $$ \det_{0\leq i,j\leq n} \left( {E/2\choose i+j+1}
        \right). $$ 
This last determinant is well-known; see e.g.\ \cite[(3.12)]{kratt}. 
\qed

Stanton has pointed out that the determinant of (\ref{eq:oddps}) is a
special case of a Hankel determinant of Meixner polynomials
$M_n(x;b,c)$, viz., $a_p(E)=2EM_{p-1}(E-1;2,-1)$. Since the Meixner
polynomials are moments of a Jacobi polynomial measure
\cite[Thm.~524]{i-s} the determinant will explicitly factor. 

Neither of these two proofs of factorization is very enlightening. Is
there a more conceptual proof based on the theory of symmetric
functions?

\textsc{Note.} Permutations whose shape is a $p\times p$ square have
an alternative description as a consequence of a basic property of the
RSK algorithm \cite[Cor.~7.23.11, Thm.~7.23.17]{ec2}, viz., they are
the permutations in $\fs_{p^2}$ whose longest increasing subsequence
and longest decreasing subsequence both have length $p$.

There are some other ``special factorizations'' of $\alt(\lambda)$ and
$\ralt(\lambda)$ that appear to hold, which undoubtedly can be proved
in a manner similar to the proof of Theorem~\ref{thm:opsdet}. Some of
these cases are the following, together with those arising from the
identity $\alt(\lambda)= \alt(\lambda')$ when $|\lambda|$ is odd, and
$\alt(\lambda)=\ralt(\lambda')$ when $|\lambda|$ is even. We write
$\lambda=\langle 1^{m_1}2^{m_2}\cdots\rangle$ to indicate that
$\lambda$ has $m_i$ parts equal to $i$.
 \begin{itemize}
   \item $\ralt(\langle p^{p-1}\rangle)$
   \item $\alt(\langle 1,p^p\rangle)$, $p$ odd
   \item certain values of $\alt(b,b-1,b-2,\dots,a)$ or
   $\ralt(b,b-1,b-2,\dots,a)$.
 \end{itemize}
There are numerous other values of $\lambda$ for which $\alt(\lambda)$
or $\ralt(\lambda)$ ``partially factors.'' Moreover, there are similar
specializations of $s_\lambda$ which factor nicely, although they
don't correspond to values of $\alt(\lambda)$ or $\ralt(\lambda)$,
e.g., $s_{\langle p^p\rangle}[E,0,-E,0,E,0,-E,0,\dots]$ for $p$ even. 

\section{Cycle type.}
A permutation $w\in\sn$ has \emph{cycle type} $\rho(w)=(\rho_1,\rho_2,
\dots)\vdash n$ if the cycle lengths of $w$ are $\rho_1,
\rho_2,\dots$. For instance, the identity permutation has cycle type
$\langle 1^n\rangle$. In this section we give an umbral formula for
the number of alternating and reverse alternating permutations
$w\in\sn$ of a fixed cycle type. 

Our results are based on a theorem of Gessel-Reutenauer \cite{g-r},
which we now explain. Define a symmetric function
  \beq L_n = \frac 1n\sum_{d\mid n}\mu(d)p_d^{n/d}, \label{eq:lndef}
   \eeq
where $\mu$ is the number-theoretic M\"obius function.
Next define $L_{\langle m^r\rangle}=h_r[L_m]$
(plethysm). Equivalently, if $f(x)=f(x_1,x_2,\cdots)$ then write
$f(x^r) =f(x_1^r,x_2^r,\cdots)$. Then for fixed $m$ we have
  \beq \sum_{r\geq 0} L_{\langle m^r\rangle}(x)t^r
    = \exp \sum_{r\geq 1} \frac 1r L_m(x^r)t^r. 
    \label{eq:lmr} \eeq
Finally, for any partition $\lambda=\langle 1^{m_1}2^{m_2}\cdots
\rangle$ set 
  \beq L_\lambda=L_{\langle 1^{m_1}\rangle} 
         L_{\langle 2^{m_2}\rangle}\cdots. \label{eq:mult} \eeq 
For some properties of the symmetric functions $L_\lambda$ see
\cite[Exer.~7.89]{ec2}.

\begin{theorem}[Gessel-Reutenauer] \label{thm:gr}
Let $\rho\vdash n$ and $\alpha\in\con$. Let $f(\rho,\alpha)$ denote
the number of permutations $w\in\sn$ satisfying $\rho=\rho(w)$ and
$\alpha= \co(w)$. Then 
  $$ f(\rho,\alpha) = \langle L_\rho, s_{B_\alpha}\rangle. $$
\end{theorem}

Now for $\rho\vdash n$ let $b(\rho)$ (respectively, $b^\ast(\rho))$
denote the number of alternating (respectively, reverse alternating)
permutations $w\in\sn$ of cycle type $\rho$. The following corollary
is then the special cases $B_\alpha=\tau_n$ and $B_\alpha=\tau'_n$ of
Theorem~\ref{thm:gr}.

\begin{corollary} \label{cor:brho}
We have $b(\rho) =\langle L_\rho, s_{\tau_n}\rangle$ and $b^\ast(\rho)
=\langle L_\rho, s_{\tau'_n}\rangle$. 
\end{corollary}

We first consider the case when $\rho=(n)$, i.e., $w$ is an
$n$-cycle. Write $b(n)$ and $b^\ast(n)$ as short for $b((n))$ and
$b^\ast((n))$. Theorem~\ref{thm:bn} below is actually subsumed by 
subsequent results (Theorems~\ref{thm:moddfm} and \ref{thm:mevenfm}),
but it seems worthwhile to state it separately.

\begin{theorem} \label{thm:bn}
\emph{(a)} If $n$ is odd then
  $$ b(n) = b^\ast(n) = \frac 1n \sum_{d\mid
    n}\mu(d)(-1)^{(d-1)/2}E_{n/d}. $$ 
\emph{(b)} If $n=2^km$ where $k\geq 1$, $m$ is odd, and $m\geq 3$,
then
  $$ b(n) = b^*(n) = \frac 1n \sum_{d\mid m}\mu(d)E_{n/d}. $$
\emph{(c)} If $n=2^k$ and $k\geq 2$ then
  \beq b(n) = b^\ast(n) = \frac 1n(E_n-1). \label{eq:bnc} \eeq
\emph{(d)} Finally, $b(2)=1$, $b^\ast(2)=0$. 
\end{theorem}

\proof (a) By Theorem~\ref{thm:main} and Corollary~\ref{cor:brho} we
have for odd $n$ that 
  \beas b(n) = b^\ast(n) & = & L_n[E,0,-E,0,E,0,-E,0,\cdots]\\
   & = & \frac 1n\sum_{d\mid n} \mu(d)((-1)^{(d-1)/2}E)^{n/d}\\
   & = & \frac 1n \sum_{d\mid n}\mu(d)(-1)^{(d-1)/2}E_{n/d}, \eeas
since $n/d$ is odd for each $d\mid n$.

(b) Split the sum (\ref{eq:lndef}) into two parts: $d$ odd and $d$
even. Since $\mu(2d)=-\mu(d)$ when $d$ is odd and since $\mu(4d)=0$
for any $d$,  we obtain
  \beas b(n) & = & L_n[E,-1,-E,1,E,-1,-E,1,\cdots]\\ & = &
    \frac 1n\left(\sum_{d\mid m} \mu(d)((-1)^{(d-1)/2}E)^{n/d} -
      \sum_{d\mid m}\mu(d)((-1)^d)^{n/2d}\right)\\ & = &
     \frac 1n\left(\sum_{d\mid m} \mu(d)E_{n/d} -
      (-1)^{n/2}\sum_{d\mid m}\mu(d)\right). \eeas
The latter sum is 0 since $m>1$, and we obtain the desired formula for 
$b(n)$. The argument for $b^\ast(n)$ is completely analogous; the factor
$(-1)^{n/2}$ now becomes $(-1)^{1+\frac n2}$.

(c) When $n=2^k$, $k\geq 2  $,  we have
  $$ L_n =\frac 1n\left( p_1^n-p_2^{n/2}\right). $$
Substituting $p_1=E$ and $p_2=\pm 1$, and using that $n/2$ is even,
yields (\ref{eq:bnc}).

(d) Trivial. It is curious that only for $n=2$ do we have $b(n)\neq
b^\ast(n)$. \qed

Note the special case of Theorem~\ref{thm:bn}(a) when $m=p^k$, where
$p$ is an odd prime and $k\geq 1$:
  $$ b(p^k) =\frac{1}{p^k}\left(E_{p^k}-(-1)^{(p-1)/2}
       E_{p^{k-1}}\right). $$ 
Is there a simple combinatorial proof, at least when $k=1$? The same
can be asked of equation~(\ref{eq:bnc}).
  
We next turn to the case $\lambda=\langle m^r\rangle$, i.e., all
cycles of $w$ have length $m$. Write $b(m^r)$ as short for $b(\langle
m^r\rangle)$, and similarly for $b^\ast(m^r)$. Set
  \beas F_m(t) & = & \sum_{r\geq 0}b(m^r)t^r\\
       F^\ast_m(t) & = & \sum_{r\geq 0}b^\ast(m^r)t^r. \eeas
First we consider the case when $m$ is odd.

\begin{theorem} \label{thm:moddfm}
\emph{(a)} Let $m$ be odd and $m\geq 3$. Then
  $$ F_m(t) = F^\ast_m(t) = \exp\left[ \frac 1m
   \left(\sum_{d\mid m}\mu(d)(-1)^{(d-1)/2}E^{m/d}\right)
   (\tan^{-1}t)\right]. $$
\emph{(b)} We have
  \beas F_1(t) & = & \sinh(E\tan^{-1} t) +
   \frac{1}{\sqrt{1+t^2}}\cosh(E\tan^{-1}t)\\[.05in]
   F_1^\ast(t) & = & \sinh(E\tan^{-1} t) +
   \sqrt{1+t^2}\cosh(E\tan^{-1}t).
  \eeas
\end{theorem}

\proof (a) By equations~(\ref{eq:feo}) and (\ref{eq:lmr}) we have that
the terms of $F_m(t)$ and $F^\ast_m(t)$ of odd degree (in $t$) are
given by 
  \bea \frac 12(F_m(t) -F_m(-t)) & = & \frac 12(F_m^\ast(t)-
    F_m^\ast(-t))\nonumber\\ & = &
  \left(\sinh \sum_{r\odd} \frac 1r   
    L(x^r)t^r\right)\nonumber \\ 
   & & \ \ \left( \exp\sum_{r\even} \frac 1r    
    L(x^r)t^r\right)[E,0,-E,0,\dots]\nonumber\\ & = &
   \left( \sinh \sum_{r\odd}\frac{t^r}{mr}\sum_{d\mid m}
    \mu(d)p_{rd}^{m/d}\right)[E,0,-E,0,\dots]\\ & = &
    \sinh \sum_{r\odd} \frac{t^r}{mr}\sum_{d\mid m}
    \mu(d)(-1)^{(rd-1)/2}E^{m/d}\nonumber\\ &  = & 
    \sinh \frac 1m\sum_{d\mid m} \mu(d)(-1)^{(d-1)/2}E^{m/d}
       \left(t-\frac{t^3}{3}+\frac{t^5}{5}-\cdots\right)
    \nonumber\\ & = &
    \sinh \frac 1m\left(\sum_{d\mid m} \mu(d)(-1)^{(d-1)/2}
     E^{m/d}\right) (\tan^{-1} t). \label{eq:sinh} \eea

Similarly the terms of $F_m(t)$ of even degree are given by
  \beas \frac 12(F_m(t) +F_m(-t))  & = &
  \left(\cosh \sum_{r\odd} \frac 1r   
    L(x^r)t^r\right)\\ 
    & & \ \ \cdot\left( \exp\sum_{r\even} \frac 1r   
    L(x^r)t^r\right)[E,-1,-E,1,\dots]\nonumber\\ & = &
   \left( \cosh \sum_{r\odd}\frac{t^r}{mr}\sum_{d\mid m}
    \mu(d)p_{rd}^{m/d}\right)\\ 
     & & \ \ \cdot\left( \exp \sum_{r\even}
     \frac{t^r}{mr}\sum_{d\mid m}\mu(d)p_{rd}^{m/d}\right)
    [E,-1,-E,1,\dots]\\ & = &
    \left(\cosh\sum_{r\,\mathrm{odd}} \frac{t^r}{mr}\sum_{d\mid m}
    \mu(d)((-1)^{(rd-1)/2})^{m/d}E^{m/d}\right) \\ 
      & & \ \ \cdot\left(\exp\sum_{r\,\mathrm{even}}\frac{t^r} 
     {mr}\sum_{d\mid m}\mu(d)((-1)^{rd/2})^{m/d}\right)
    \\ & = &
    \left(\cosh \frac 1m\sum_{d\mid m} \mu(d)(-1)^{(d-1)/2}E^{m/d}
       \tan^{-1}t\right) \\ & & \ \
     \cdot\left(\exp\sum_{r\,\mathrm{even}}\frac{t^r}  
     {mr}(-1)^{r/2}\sum_{d\mid m}\mu(d) \right) \eeas
    \beq \quad = \cosh \frac 1m\left(\sum_{d\mid m}
  \mu(d)(-1)^{(d-1)/2} E^{m/d}\right)(\tan^{-1} t). \label{eq:cosh}
  \eeq
Adding equations (\ref{eq:sinh}) and (\ref{eq:cosh}) yields (a) for
$F_m(t)$. 

The computation for $F^\ast_m(x)$ is identical, except that the factor
$(-1)^{rd/2}$ is replaced by $(-1)^{1+rd/2}$. This alteration does not
affect the final answer.

(b) The computation of the odd part of $F_1(t)$ and $F^\ast_1(t)$ is
the same as in (a), yielding
  \beas \frac 12(F_1(t) -F_1(-t)) & = & \frac 12(F_1^\ast(t)-
    F_1^\ast(-t))\\ & = & \sinh (E\tan^{-1}t). \eeas
On the other hand,
  \beas \frac 12(F_1(t) +F_1(-t))  & = &  \cosh( E\tan^{-1}t)
    \cdot\left( \exp\sum_{r\even}\frac{t^r}{r}\mu(1)(-1)^{r/2}
    \right)\\ & = & \frac{\cosh(E\tan^{-1}t)}
    {\sqrt{1+t^2}}, \eeas
and the proof for $F_1(t)$ follows. For $F_1^\ast(t)$ the factor 
$(-1)^{r/2}$ becomes $(-1)^{1+r/2}$, so the factor $\sqrt{1+t^2}$
moves from the denominator to the numerator.
\qed

Clearly the only alternating permutation of cycle type $\langle
1^r\rangle$ is $1$ (when $r=1$). Hence from
Theorem~\ref{thm:moddfm}(b) we obtain the umbral identity
  \beq \sinh(E\tan^{-1}t)+\frac{1}{\sqrt{1+t^2}}\cosh(E\tan^{-1}t)
    = 1+t. \label{eq:umbralid} \eeq
One may wonder what is the point of Theorem~\ref{thm:moddfm}(b) since
$b(1^r)$ is trivial to compute directly. Its usefulness will be seen
below (Theorem~\ref{thm:cycind}), when we consider ``mixed'' cycle
types, i.e., not all cycle lengths are equal.

Theorem~\ref{thm:moddfm} can be restated ``non-umbrally'' analogously
to Theorem~\ref{thm:doubalt}. For instance, if $m=p^k$ where $p$ is
prime and $p\equiv 3\,(\mathrm{mod}\,4)$, then
  $$ F_m(t) = F^\ast_m(t) = \sum_{i,j\geq 0} E_{(m/p)i+mj}
    \frac{\left( \frac 1p\tan^{-1}t\right)^{i+j}}{i!\,j!}, $$
while if $p\equiv 1\,(\mathrm{mod}\,4)$, then
  $$ F_m(t) = F^\ast_m(t) = \sum_{i,j\geq 0} (-1)^iE_{(m/p)i+mj}
    \frac{\left( \frac 1p\tan^{-1}t\right)^{i+j}}{i!\,j!}. $$
For general odd $m$, $F_m(t)$ will be expressed as a $2^{\nu(m)}$-fold 
sum, where $\nu(m)$ is the number of distinct prime divisors of $m$.

\begin{theorem} \label{thm:mevenfm}
\emph{(a)} Let $m=2^kh$, where $k\geq 1$, $h\geq 3$, and $h$ is
odd. Then 
  $$ F_m(t)=F_m^\ast(t) =\left( \frac{1+t}{1-t}\right)
      ^{\frac{1}{2m}\sum_{d\mid h}\mu(d)E^{m/d}}. $$
\emph{(b)} Let $m=2^k$ where $k\geq 2$. Then
  $$ F_m(t) = F_m^\ast(t) = \left(\frac{1+t}{1-t}\right)
    ^{\frac{1}{2m}(E^m-1)}. $$
 \emph{(c)} Let $m=2$. Then
  \beas F_2(t) & = &  \left(\frac{1+t}{1-t}
    \right)^{(E^2+1)/4}\\[.05in] F^\ast_2(t) & = & \displaystyle
    \frac{F_2(t)}{1+t}\ \ (\mathrm{compare}\ (\ref{eq:aa4})). \eeas
\end{theorem}

\proof (a) The argument is analogous to the proof of
Theorem~\ref{thm:moddfm}. We have
  \beas F_m(t) & = & \left(\exp \sum_{r\geq 1} \frac 1r  
    L(x^r)t^r\right)[E,-1,-E,1,\dots]\\ & = &
   \left( \exp \sum_{r\geq 1}\frac{t^r}{mr}\sum_{d\mid m}
    \mu(d)p_{rd}^{m/d}\right)[E,-1,-E,1,\dots]\\ & = &
   \exp\left( \sum_{r\,\mathrm{odd}}\frac{t^r}{rm}
   \sum_{d\mid h} ((-1)^{(rd-1)/2})^{m/d}\mu(d)E^{m/d}\right.\\ 
    & & \ + \left.\sum_{r\,\mathrm{even}}\frac{t^r}{rm}
   \sum_{d\mid h} ((-1)^{rd/2})^{m/d} -\sum_r\frac{t^r}{rm}
   \sum_{d\mid h} ((-1)^{rd/2})^{m/d}\right)\\ & = &
   \exp \sum_{r\,\mathrm{odd}}\frac{t^r}{rm}\sum_{d\mid h}
   \mu(d)E^{m/d}\\ & = & \exp \left(\frac 1m\sum_{d\mid h}
   \mu(d)E^{m/d}\right)\frac 12\log\frac{1+t}{1-t}, \eeas
and the proof follows for $F_m(t)$. The same argument holds for
$F^\ast_m(t)$ since $-1$ was always raised to an even power or was 
multiplied by a factor $\sum_{d\mid h}\mu(d)=0$ in the proof. 

(b) We now have 
  \beas F_m(t) & = & \exp \sum_{r\geq 1}\frac{t^r}{rm}\left(
p_r^m-p_{2r}^{m/2}\right)[E,-1,-E,1,\dots]\\ & = &
  \exp \frac 1m\left(\sum_{r\,\mathrm{odd}}\frac{t^r}{r}
   \left( ((-1)^{(r-1)/2}E)^m-(-1)^{rm/2}\right)\right.\\ & & \ + 
    \left.\sum_{r\,\mathrm{even}}\frac{t^r}{r}\left((-1)^{rm/2}
     -(-1)^{rm/2}\right)\right)\\ & = &
    \left(\frac{1-t}{1+t}\right)^{1/2m}\exp \frac{E^m}{2m}
    \log\frac{1+t}{1-t}, \eeas
etc. Again the computation for $F^\ast_m(t)$ is the same.

(c) We have
  \beas F_2(t) & = & \exp \frac 12\sum_{r\geq 1}\left(p_r^2-p_{2r} 
    \right)\frac{t^r}{r}[E,-1,-E,1,\dots]\\ & = &
  \exp \frac 12\left[ \sum_{r\,\mathrm{odd}}\left(
   \left( (-1)^{(r-1)/2}E\right)^2-(-1)^r\right)
  \frac{t^r}{r}\right.\\ & & \ +\left. \sum_{r\,\mathrm{even}}
   \left((-1)^r-(-1)^r\right)\right]\\
   & = & \exp \frac 12\sum_{r\,\mathrm{odd}}(E^2+1)\frac{t^r}{r}.
   \eeas
etc. We leave the case $F_2^\ast(t)$ to the reader. \qed

The expansion of $F_2(t)$ begins
  \beq F_2(t) = 1+t+t^2+2t^3+5t^4+17t^5+72t^6+367t^7+2179t^8+\cdots. 
    \label{eq:ram} \eeq
Ramanujan asserts in Entry 16 of his second notebook (see
\cite[p.~545]{berndt}) that as $t$ tends to $0+$,
  \beq 2\sum_{n\geq 0}(-1)^n\left( \frac{1-t}{1+t}\right)^{n(n+1)}
    \sim 1+t+t^2+2t^3+5t^4+17t^5+\cdots. \label{eq:berndt} \eeq
Berndt \cite[(16.6)]{berndt} obtains a formula for the complete
asymptotic expansion of $2\sum_{n\geq 0}(-1)^n\left(
\frac{1-t}{1+t}\right)^{n(n+1)}$ as $t\rightarrow 0+$. It is easy to
see that Berndt's formula can be written as $\left(\frac{1+t}{1-t}
\right)^{(E^2+1)/4}$ and is thus equal to
$F_2(t)$. Theorem~\ref{thm:mevenfm}(c) therefore answers a question of
Galway \cite[p.~111]{galway}, who asks for a combinatorial
interpretation of the coefficients in Ramanujan's asymptotic
expansion. 

\textsc{Note.} The following formula for $F_2(t)$ follows from
equation (\ref{eq:berndt}) and an identity of Ramanujan proved by
Andrews \cite[(6.3)$_\mathrm{R}$]{andrews}: 
  $$ F_2(t) = 2\sum_{n\geq 0}q^n\frac{\prod_{j=1}^n(1-q^{2j-1})}
    {\prod_{j=1}^{2n+1}(1+q^j)}, $$
where $q=\left(\frac{1-t}{1+t}\right)^{2/3}$. It is not hard to see
that this is a \emph{formal} identity, unlike the asymptotic identity
(\ref{eq:berndt}). 

\textsc{Note.} We can put Theorems~\ref{thm:moddfm}(a) into a form
more similar to Theorem~\ref{thm:mevenfm}(a) by noting the identity
  $$ \exp (\tan^{-1} t) = \left( \frac{1-it}{1+it}\right)^{i/2}. $$
Hence when $m$ is odd and $m\geq 3$ we have
  $$ F_m(t)=F^\ast_m(t) =\left( \frac{1-it}{1+it}\right)^
    {\frac{i}{2m}
   \sum_{d\mid m}\mu(d)(-1)^{(d-1)/2}E^{m/d}}. $$

The multiplicativity property (\ref{eq:mult}) of $L_\lambda$ allows us
write down a generating function for the number $b(\lambda)$
(respectively, $b^\ast(\lambda)$) of alternating (respectively,
reverse alternating) permutations of any cycle type $\lambda$. For
this purpose, let $t_1,t_2,\dots$ and $t$ be indeterminates and set
$\deg(t_i)=i$, $\deg(t)=1$. If $F(t_1,t_2,\dots)$ is a power series in 
$t_1,t_2,\dots$ or $F(t)$ is a power series in $t$, then write 
${\cal O}F$ (respectively, ${\cal E}F$) for those terms of $F$ whose
total degree is odd (respectively, even). For instance,
  $$ {\cal O}F(t_1,t_2,\dots) = \frac 12\left(F(t_1,t_2,t_3,t_4\dots)-
      F(-t_1,t_2,-t_3,t_4,\dots)\right). $$
Define the ``cycle indicators''
  \beas Z(t_1,t_2,\dots) & = & \sum_{\lambda=\langle 1^{m_1}2^{m_2} 
    \cdots\rangle}  b(\lambda) t_1^{m_1} t_2^{m_2}\cdots\\
    Z^\ast(t_1,t_2,\dots) & = & \sum_{\lambda=\langle 1^{m_1}2^{m_2} 
    \cdots\rangle}  b^\ast(\lambda) t_1^{m_1} t_2^{m_2}\cdots, \eeas
where both sums range over all partitions $\lambda$ of all integers
$n\geq 0$.

\begin{theorem} \label{thm:cycind}
We have
  $$ \hspace{-2in}{\cal O}Z(t_1,t_2,\dots) = 
    {\cal O}Z^\ast(t_1,t_2,\dots) $$
  \beq =
     {\cal O}\exp
  (E\tan^{-1}t_1)\cdot\left(\frac{1+t_2}{1-t_2}\right)^{E^2/4} 
     F_3(t_3)F_4(t_4)\cdots \label{eq:zodd} \eeq
   \beas {\cal E}Z(t_1,t_2,\dots) & = & {\cal E}
     \frac{\exp(E\tan^{-1}t_1)}{\sqrt{1+t_1^2}}
    \left(\frac{1+t_2}{1-t_2}\right)^{(E^2+1)/4} 
     F_3(t_3)F_4(t_4)\cdots \\[.1in]
   {\cal E}Z^\ast(t_1,t_2,\dots) & = & {\cal E}
     \sqrt{1+t_1^2}\,\exp(E\tan^{-1}t_1)\cdot
    \frac{1}{1+t_2}\left(\frac{1+t_2}{1-t_2}\right)^{(E^2+1)/4}\\ 
       & & \ \ \cdot F_3(t_3)F_4(t_4)\cdots. \eeas
It is understood that in these formulas $F_j(t_j)$ is to be written
in the umbral form given by Theorems~\ref{thm:moddfm} and
\ref{thm:mevenfm}. 
\end{theorem}

\proof Let 
  \beq G_m(t) = \exp \sum_{r\geq 1} \frac 1r L_m(x^r)t^r. 
   \label{eq:gmt} \eeq
It follows from equations (\ref{eq:lmr}) and (\ref{eq:mult})
that 
  $$ {\cal O}Z(t_1,t_2,\dots) = {\cal O}\prod_{m\geq 1}
     G_m(t_m)[E,0,-E,0,\dots]. $$
The proofs of Theorems~\ref{thm:moddfm}(a) and \ref{thm:mevenfm}(a,b)
show that for $m\geq 3$, 
  $$ G_m(t)[E,0,-E,0,\dots] = G_m(t)[E,-1,-E,1,\dots]. $$
Hence we obtain the factors $F_3(t_3)F_4(t_4)\cdots$ in 
equation~(\ref{eq:zodd}). It is straightforward to compute 
$G_m(t_m)[E,0,-E,0,\dots]$ (and is implicit in
the proofs of Theorems~\ref{thm:moddfm}(b) and \ref{thm:mevenfm}(c))
for $m=1,2$. For instance, 
  \bea G_1(t_1)[E,0,-E,0\dots] & = & \exp\left(\sum_{r\geq 1} \frac 1r
   p_rt_1^r\right)[E,0,-E,0,\dots] \nonumber\\ & = &
    \exp \sum_{r\odd}\frac 1r(t_1-\frac 13t_1^3+\cdots) \nonumber\\
     & = & \exp(\tan^{-1} t_1). \label{eq:eatt} \eea
Thus we obtain the remaining factors in equation~(\ref{eq:zodd}). The
remaining formulas are proved analogously. \qed

We mentioned in Section~\ref{sec1} that Ehrenborg and Readdy raised
the question of counting alternating involutions $w\in S_n$. An answer
to this question is a simple consequence of Theorem~\ref{thm:cycind}. 

\begin{corollary}
Let $c(n)$ (repectively, $c^\ast(n)$) denote the number of alternating
(respectively, reverse alternating) involutions $w\in\sn$. Then  
  \beas \sum_{n\geq 0} c(2n+1)t^{2n+1} & = & \sinh(E\tan^{-1}t)\cdot 
           \left(\frac{1+t^2}{1-t^2}\right)^{E^2/4}\\[.05in] 
    \sum_{n\geq 0} c(2n)t^{2n} & = & \frac{1}{\sqrt[4]{1-t^4}}\,
      \cosh(E\tan^{-1}t)\cdot
           \left(\frac{1+t^2}{1-t^2}\right)^{E^2/4}\\[.05in] 
    c^\ast(n) & = & c(n). \eeas    
Equivalently, 
  \beas \sum_{n\geq 0} c(2n+1)t^{2n+1} & = &   
 \sum_{i,j\geq 0} \frac{E_{2i+2j+1}}{(2i+1)!\,j!\,4^j}
   \tan^{-1}(t)^{2i+1} \left(\log \frac{1+t^2}
   {1-t^2}\right)^j \\[.05in]
   \sum_{n\geq 0} c(2n)t^{2n} & = & \frac{1}{\sqrt[4]{1-t^4}}\,   
    \sum_{i,j\geq 0} \frac{E_{2i+2j}}{(2i)!\,j!\,4^j}
    \tan^{-1}(t)^{2i} \left(\log \frac{1+t^2}
   {1-t^2}\right)^j. \eeas
\end{corollary}

\proof
We have
  \beas \sum_{n\geq 0}c(2n+1)t^{2n+1} & = & {\cal O}
      Z(t,t^2,0,0,\dots) \\
   \sum_{n\geq 0}c(2n)t^{2n} & = & {\cal E}
      Z(t,t^2,0,0,\dots), \eeas
and similarly for $c^\ast(n)$. The result is thus a special case of
Theorem~\ref{thm:cycind}. \qed

The identity $c(n)=c^\ast(n)$ does not seem obvious. It can also be
obtained using properties of the RSK algorithm, analogous to the
argument after the proof of Theorem~\ref{thm:doubalt} Namely,
$w\in\sn$ is an alternating (respectively, reverse alternating)
involution if and only if $w\rsk(P,P)$, where $P$ is an alternating
(respectively, reverse alternating) SYT. Hence if $w'\rsk (P^t,P^t)$,
then the map $w\mapsto w'$ interchanges alternating involutions
$w\in\sn$ with reverse alternating involutions $w'\in\sn$.

\section{Fixed points.}
P. Diaconis (private communication) raised the question of enumerating
alternating permutations by their number of fixed points. It is easy
to answer this question using Theorem~\ref{thm:cycind}. Write $d_k(n)$
(respectively, $d^*_k(n)$) for the number of alternating
(respectively, reverse alternating) permutations in $\sn$ with $k$
fixed points.  Write ${\cal O}_t$ and ${\cal E}_t$ for the odd and
even part of a power series with respect to $t$ (ignoring other
variables), i.e., ${\cal O}_tF(t)=\frac 12(F(t)-F(-t))$ and ${\cal
E}_tF(t)=\frac 12 (F(t)+F(-t))$.

\begin{proposition} \label{prop:fixed}
We have 
  \bea \sum_{k,n\geq 0} d_k(2n+1)q^k t^{2n+1} & = & {\cal O}_t
   \frac{\exp(E(\tan^{-1}qt-\tan^{-1}t))}{1-Et} \label{eq:dkodd}\\
  d^*_k(2n+1) & = & d_k(2n+1) \label{eq:dksodd} \\
  \sum_{k,n\geq 0} d_k(2n)q^k t^{2n} & = & {\cal E}_t
  \sqrt{\frac{1+t^2}{1+q^2t^2}}
  \frac{\exp(E(\tan^{-1}qt-\tan^{-1}t))}{1-Et} \nonumber \\ 
   \sum_{k,n\geq 0} d_k^*(2n)q^k t^{2n} & = &  {\cal E}_t
  \sqrt{\frac{1+q^2t^2}{1+t^2}}
  \frac{\exp(E(\tan^{-1}qt-\tan^{-1}t))}{1-Et} \nonumber.
  \eea
Equivalently, we have the nonumbral formulas
  \beas \sum_{k,n\geq 0} d_k(2n+1)q^k t^{2n+1} & = & 
    \sum_{\twoline{i,j\geq 0}{i\not\equiv j\,(\mathrm{mod}\,2)}}
     \frac{E_{i+j}}{j!}t^i(\tan^{-1}qt-\tan^{-1}t)^j\\
   \sum_{k,n\geq 0} d_k(2n)q^k t^{2n} & = & 
    \sqrt{\frac{1+t^2}{1+q^2t^2}}
    \sum_{\twoline{i,j\geq 0}{i\equiv j\,(\mathrm{mod}\,2)}}
     \frac{E_{i+j}}{j!}t^i(\tan^{-1}qt-\tan^{-1}t)^j\\
    \sum_{k,n\geq 0} d_k^*(2n)q^k t^{2n} & = & 
    \sqrt{\frac{1+q^2t^2}{1+t^2}}
    \sum_{\twoline{i,j\geq 0}{i\equiv j\,(\mathrm{mod}\,2)}}
     \frac{E_{i+j}}{j!}t^i(\tan^{-1}qt-\tan^{-1}t)^j.
  \eeas
\end{proposition}

\proof 
It is not hard to see (e.g., \cite[(1)]{schocker}) that
  $$ \sum_{\lambda\vdash n} L_\lambda = p_1^n, $$
where $p_1=x_1+x_2+\cdots$. It follows from equations ({\ref{eq:lmr}),
(\ref{eq:mult}) and (\ref{eq:gmt}) that 
  $$ G_1(t)G_2(t)\cdots = \sum_{n\geq 0} p_1^nt^n =
  \frac{1}{1-p_1t}. $$ 
Hence by equation (\ref{eq:zodd}) we have
  \beas \sum_{k,n\geq 0}d_k(2n+1)q^kt^{2n+1} & = &
    \sum_{k,n\geq 0}d_k^*(2n+1)q^kt^{2n+1}\\ & = &
   {\cal O}_t \exp(E \at qt)\left(\frac{1+t}{1-t}\right)^{E^2/4} 
     F_3(t)F_4(t)\cdots\\ & = &
   {\cal O}_t \frac{\exp(E \at qt)}{\exp(E\at t)\cdot(1-Et)}, \eeas
proving (\ref{eq:dkodd}) and
  (\ref{eq:dksodd}). The proof for $n$ even is analogous.
\qed

\begin{corollary}
For $n>1$ we have $d_0(n)=d_1(n)$ and $d_0^*(n)=d_1^*(n)$.
\end{corollary}

\proof Let 
  $$ M(q,t) = {\cal O}_t \frac{\exp
       E(\tan^{-1}qt-\tan^{-1}t)}{1-Et}. $$
By equation~(\ref{eq:dksodd}) it follows that
  \beas \sum_{n\,\mathrm{odd}} d_0(n)t^n & = &
       M(0,t)\\ 
    \sum_{n\,\mathrm{odd}} d_1(n)t^n & = &
       \left.\frac{\partial}{\partial q} M(q,t)\right|_{q=0}. \eeas
It is straightforward to compute that
   $$ \left.\frac{\partial}{\partial q}
        M(q,t)\right|_{q=0}-M(0,t) = \sinh(E\tan^{-1} t). $$
By equation~(\ref{eq:umbralid}) we have $\sinh(E\tan^{-1} t)=t$, and 
the proof follows for $n$ odd. The proof for $n$ even is completely
analogous.  
\qed

We have a conjecture about certain values of $d_k(n)$ and
$d_k^*(n)$. It is not hard to see that 
  \beas \max\{k\st d_k(n)\neq 0\} & = & \lceil n/2\rceil,\quad n\geq
      4\\ 
     \max\{k\st d_k^*(n)\neq 0\} & = & \lceil (n+1)/2\rceil,\quad
     n\geq 5.  
  \eeas
 
\begin{conjecture}
Let $D_n$ denote the number of derangements (permutations
without fixed points) in $\sn$. Then
  \beas d_{\lceil n/2\rceil}(n) & = & D_{\lfloor n/2\rfloor},
   \quad n\geq 4\\ 
   d_{\lceil (n+1)/2\rceil}^*(n) & = & D_{\lfloor (n-1)/2\rfloor}, 
   \quad n\geq 5. \eeas
\end{conjecture}

It is also possible to obtain asymptotic information from
Proposition~\ref{prop:fixed}. The next result considers alternating
or reverse alternating derangements (permutations without fixed
points). 

\begin{corollary} \label{cor:asy}
\rm{(a)} We have for $n$ odd the asymptotic expansion
  \bea d_0(n) & \sim & \frac 1e\left( E_n + a_1 E_{n-2} + a_2
  E_{n-4}+\cdots\right) \label{eq:asymodd}\\ & = & \frac 1e\left(
  E_n+\frac 13 E_{n-2}- 
   \frac{13}{90}E_{n-4}+\frac{467}{5760}E_{n-6}+\cdots\right), 
    \nonumber \eea
where
   $$ \sum_{k\geq 0}a_k x^{2k} = \exp\left( 1-\frac 1x \tan^{-1}x
      \right). $$

\rm{(b)} We have for $n$ even the asymptotic expansion
 \bea d_0(n) & \sim & \frac 1e\left( E_n + b_1 E_{n-2} + b_2
  E_{n-4}+\cdots\right) \label{eq:asymeven}\\ & = & \frac 1e\left(
  E_n+\frac 56 E_{n-2}- 
   \frac{37}{360}E_{n-4}+\frac{281}{9072}E_{n-6}+\cdots\right), 
    \nonumber \eea
where
   $$ \sum_{k\geq 0}b_k x^{2k} = \sqrt{1+x^2}\exp\left( 1-\frac 1x
      \tan^{-1}x \right). $$

\rm{(c)} We have for $n$ even the asymptotic expansion
 \bea d_0^*(n) & \sim & \frac 1e\left( E_n + c_1 E_{n-2} + c_2
  E_{n-4}+\cdots\right) \label{eq:asymevens}\\ & = & \frac 1e\left(
  E_n-\frac 16 E_{n-2}+ 
   \frac{23}{360}E_{n-4}-\frac{1493}{45360}E_{n-6}+\cdots\right), 
    \nonumber \eea
where
   $$ \sum_{k\geq 0}c_k x^{2k+1} = \frac{1}{\sqrt{1+x^2}}\exp\left(
      1-\frac 1x \tan^{-1}x \right). $$
\end{corollary}

\textsc{Note.}  Equations~(\ref{eq:asymodd}), (\ref{eq:asymeven}), and
(\ref{eq:asymevens}) are genuine asymptotic expansions since $E_m \sim
2(2/\pi)^{m+1}m!$, 
so for fixed $k$,
  $$ E_{n-k} \sim 2\left( \frac{\pi}{2}\right)^k\frac{1}{n^k}E_n $$
as $n\rightarrow\infty$. In fact, since 
  $$ E_m = 2\left(\frac{2}{\pi}\right)^{m+1}m!(1+O(3^{-m})), $$
we can rewrite (\ref{eq:asymodd}) (and similarly (\ref{eq:asymeven}) and
(\ref{eq:asymevens})) as
  $$ d_0(n) \sim \frac{E_n}{e}\left( 1+a_1
    \left(\frac{\pi}{2}\right)^2 \frac{1}{(n)_2} +a_2
    \left(\frac{\pi}{2}\right)^4 \frac{1}{(n)_4}+\cdots\right), $$
where $(n)_j=n(n-1)\cdots(n-j+1)$. 

\emph{Proof of Corollary~\ref{cor:asy}.} (a) It follows from
equation~(\ref{eq:dkodd}) that 
  $$ \sum_{n\,\mathrm{odd}} d_0(n)t^n = {\cal O}_t
   \frac{\exp(-E\tan^{-1}t)}{1-Et}. $$
This series has the form
  $$ \sum_{n\,\mathrm{odd}} t^n(a_{n0}E^n+a_{n1} E^{n-2} + a_{n2}
  E^{n-4}+\cdots). $$
If we replace $t$ with $Et$ and $E$ with $1/E$ we therefore obtain 
  \beq {\cal O}_t \frac{\exp(-E^{-1}\tan^{-1}tE))}{1-t} =
     \sum_{n\,\mathrm{odd}} t^n(a_{n0}+a_{n1} E^2 + a_{n2}
     E^4+\cdots).  \label{eq:asymp1} \eeq
We claim that for fixed $j$ the coefficients $a_{nj}$ rapidly approach
(finite) limits as $n\rightarrow \infty$. If we expand the left-hand
side of (\ref{eq:asymp1}) as a power series in $E$, it is not hard to
see that the coefficient of $E^{2j}$ has the form $Q_j(t)/(1-t^2)$,
where $Q_j(t)$ is a polynomial in $t,e^t$ and $e^{-t}$. Hence
the coefficient of $t^{2n+1}$ in $Q_j(t)$ has the form
$p_j(n)/(2n+1)!$ for some polynomial $p_j(n)$. It follows that
   $$ a_{nj} = Q_j(1)+o(n^{-r}) $$
for all $r>0$. Now
  $$ {\cal O}_t \frac{\exp(-E^{-1}\tan^{-1}tE))}{1-t} = 
    \frac{(1+t)e^{-\frac 1E\tan^{-1}tE}-(1-t)e^{\frac 1E
     \tan^{-1} tE}}{2(1-t^2)}. $$
Multiplying by $1-t^2$ and setting $t=1$ gives $e^{-\frac 1E \tan^{-1} 
E}$, and the proof follows. The argument for (b) and (c) is
analogous. 
\qed

\section{Multisets.}
In this section we give simple umbral formulas for the number of
alternating and reverse alternating permutations of a multiset of
positive integers, with various interpretations of the meaning of
``alternating.'' There has been some previous work on alternating
multiset permutations. Goulden and Jackson \cite[Exer.~4.2.2(b);
solution, pp.~459--460]{g-j} obtain a formula for the number of
alternating permutations of the multiset with one occurrence of $i$
for $1\leq i\leq m$ and two occurrences of $i$ for $m+1\leq i\leq
m+n$. Gessel \cite[pp.~265--266]{gessel} extends this result to
multisets with one, two, or three multiplicities of each part, or with
one or four multiplicities of each part. Upon being told about the
results in this section, Gessel (private communication) was able to
extend his argument to arbitrary multisets, obtaining a result
equivalent to the case $A=\emptyset$ of Theorem~\ref{thm:multi}. Zeng
\cite{zeng} obtains an even more general result concerning the case
$A=\emptyset$. 

Our basic tool, in addition to Theorem~\ref{thm:main}, is the
following extension of Theorem~\ref{thm:dco} to skew shapes
$\lambda/\mu$. We define the descent composition of an SYT $T$ of
shape $\lambda/\mu$ exactly as for ordinary shapes, viz., $T$ has
descent composition $\alpha=(\alpha_1,\dots,\alpha_k)$ if $\{\alpha_1,
\alpha_1+\alpha_2, \dots, \alpha_1+\cdots+\alpha_{k-1}\}$ is the set
of those $i$ for which $i+1$ appears in $T$ in a lower row than $i$. 

\begin{lemma} \label{lemma:jdt}
Let $\lambda/\mu$ be a skew partition of size $n$, with corresponding
skew Schur function $s_{\lambda/\mu}$ \cite[Def.~7.10.1]{ec2}, and let
$\alpha\in\con$. Then $\langle s_{\lambda/\mu} ,s_{B_\alpha}\rangle$
is equal to the number of SYT of shape $\lambda/\mu$ and descent
composition $\alpha$.
\end{lemma}

\proof Let $s_{\lambda/\mu} =\sum_\nu c^\lambda_{\mu\nu} s_\nu$. 
Let $T$ be an SYT of shape $\lambda/\mu$, and apply jeu
de taquin \cite[{\S}A1.2]{ec2} to $T$ to obtain an SYT $T'$ of some
ordinary shape $\nu$. Two fundamental properties of jeu de
taquin assert the following:
 \begin{itemize}
  \item As $T$ runs over all SYT of shape $\lambda/\mu$,
   we obtain by jeu de taquin each SYT $T'$ of shape $\nu$
   exactly $c^\lambda_{\mu\nu}$ times. 
 \item We have $\co(T)=\co(T')$.
 \end{itemize} 
The first item above appears e.g.\ in \cite[Thm.~A1.3.1]{ec2}, while
the second item is easily proved by showing that the descent
composition is preserved by a single jeu de taquin slide. The proof
of the lemma follows immediately from the two items above. \qed

We can define $\alt(\lm)$ and $\ralt(\lm)$ for skew shapes $\lm$
exactly as we did for ordinary shapes $\lambda$. The following
corollary is then immediate from Theorem~\ref{thm:main} and
Lemma~\ref{lemma:jdt}.

\begin{corollary} \label{cor:skewalt}
Let $\lm$ be a skew shape of odd size $|\lm|$. Then
  $$ \alt(\lm) = \ralt(\lm) = s_{\lm}[E,0,-E,0,E,0,-E,0,\dots]. $$
If $|\lm|$ is even then
  \beas \alt(\lm) & = & s_{\lm}[E,-1,-E,1,E,-1,-E,1,\dots]\\
     \ralt(\lm) & = & s_{\lm}[E,1,-E,-1,E,1,-E,-1,\dots]. \eeas
\end{corollary}

We are now ready to enumerate alternating permutations of a
multiset. If two equal elements $i$ in a permutation appear
consecutively, then we need to decide whether they form an ascent or a
descent. We can make this decision separately for each $i$. Let $k\geq
1$, and let $A,B$ be complementary subsets of $[k]$, i.e.,
$A\cup B=[k]$, $A\cap B=\emptyset$. Let
$\alpha=(\alpha_1,\dots,\alpha_k)$ be a composition of some $n\geq 1$
into $k$ parts. An $\alpha$-\emph{permutation} of $[k]$ is a
permutation of the multiset
$M=\{1^{\alpha_1},\dots,k^{\alpha_k}\}$, i.e., a sequence $a_1 a_2
\cdots a_n$ with $\alpha_i$ occurrences of $i$, for $1\leq i\leq
k$. An $\alpha$-permutation is said to be $(A,B)$-\emph{alternating}
if
 $$ a_1>a_2<a_3>a_4<\cdots a_n, $$
where we define $j>j$ if $j\in A$ and $j<j$ if $j\in B$. For instance,
if $A=\{1,3\}$, $B=\{2,4\}$, and $\alpha=(3,2,2,3)$, then the
$\alpha$-permutation $w=1142214343$ is $(A,B)$-alternating since
  $$ 1>1<4>2<2>1<4>3<4>3 $$
according to our definition. Similarly we define \emph{reverse
$(A,B)$-alternating}. For example, 2213341414 is a reverse $(A,B)$
$\alpha$-permutation (with $\alpha,A,B$ as before), since
  $$ 2<2>1<3>3<4>1<4>1<4. $$
Let $N(\alpha,A,B)$ (respectively, $N^*(\alpha,A,B)$) denote the
number of $(A,B)$-alternating (respectively, reverse $(A,B)$-alternating)
$\alpha$-permutations. Write $e_i$ and $h_i$ for the elementary
and complete symmetric functions of degree $i$.

\begin{theorem} \label{thm:multi}
Let $\alpha=(\alpha_1,\dots,\alpha_k)\in\con$, and let $A,B$ be
complementary subsets of $[k]$. \\
\indent \emph{(a)}  If $n$ is odd, then 
  \beas N(\alpha,A,B) & = & N^*(\alpha,A,B)\\
   & = & \prod_{i\in A}e_{\alpha_i}\cdot\prod_{j\in B}h_{\alpha_j}
     [E,0,-E,0,E,0,-E,0,\dots]. \eeas
\indent \emph{(b)} If $n$ is even, then
  \beas N(\alpha,A,B) & = & 
   \prod_{i\in A}e_{\alpha_i}\cdot\prod_{j\in B}h_{\alpha_j}
     [E,-1,-E,1,E,-1,-E,1,\dots]\\[.1in]
      N^*(\alpha,A,B) & = & 
   \prod_{i\in A}e_{\alpha_i}\cdot\prod_{j\in B}h_{\alpha_j}
     [E,1,-E,-1,E,1,-E,-1,\dots]\\[.1in]
  \eeas
\end{theorem}

\proof Let $\sigma=\sigma(\alpha,A,B)$ be the skew
shape consisting of a disjoint union of single rows and columns, as
follows. There are $k$ connected components, of sizes
$\alpha_1,\dots,\alpha_k$ from top to bottom. If $i\in A$ then the
$i$th component is a single row, and otherwise a single column. For
instance, $\sigma((3,1,2,2),\{2,4\},\{1,3\})$ and
$\sigma((3,1,2,2),\{4\},\{1,2,3\})$ both have the following diagram:

\vspace{.5em}
\centerline{\psfig{figure=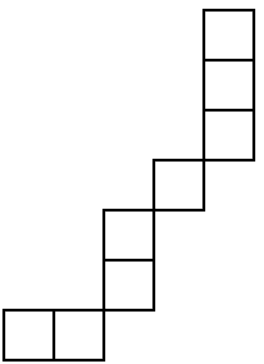}}
\vspace{.5em}

Suppose that $n$ is odd. By Corollary~\ref{cor:skewalt} we have
  $$ \alt(\sigma)=\ralt(\sigma) = s_\sigma[E,0,-E,0,\dots]. $$
Given an alternating or reverse alternationg SYT $T$ of shape
$\sigma$, define an $\alpha$-permutation $w=a_1\cdots a_n$ by the
condition that $a_i=j$ if $a_i$ appears in the $j$th component of
$\sigma$. For instance, if
  $$ \sigma = \begin{array}{cccc} & & & \!\!\!4\,8\,10\,12\\ & & 2\!\!\!\\
    & & 3\!\!\!\\ & & \!11\\ & \!5\,6\,9\\ \!1\\ \!7 \end{array}, $$
then $w=422133413121$. This construction sets up a bijection between
alternating (respectively, reverse alternating) SYT of shape $\sigma$
and $(A,B)$-alternating (respectively, reverse alternating)
$\alpha$-permutations, so the proof follows for $n$ odd. Exactly the
same argument works for $n$ even. \qed

Some values of the relevant specializations of $e_i$ and $h_i$ are as
follows:
  $$ \begin{array}{rclcl} e_1[E,0,-E,0,\dots] & = &
    h_1[E,0,-E,0,\dots] & = & E\\[.05in]
  e_2[E,0,-E,0,\dots] & = & h_2[E,0,-E,0,\dots] &
    = & \frac 12 E^2\\[.05in]
  e_3[E,0,-E,0,\dots] & = & h_3[E,0,-E,0,\dots] &
    = & \frac 16(E^3-2E)\\[.05in]
  e_4[E,0,-E,0,\dots] & = & h_4[E,0,-E,0,\dots] &
    = & \frac{1}{24}(E^4-8E^2)\\[.05in]
  e_5[E,0,-E,0,\dots] & = & h_5[E,0,-E,0,\dots] &
    = & \frac{1}{120}(E^5-20E^3+24E) \end{array} $$
 $$ \hspace{-.3in}e_1[E,-1,-E,1,\dots]=e_1[E,1,-E,-1,\dots] $$
\vspace{-.2in}
  $$ \qquad\qquad =h_1[E,-1,-E,1,\dots]=h_1[E,1,-E,-1,\dots]=E
  $$ 

  $$ \begin{array}{rclcl} e_2[E,-1,-E,1,\dots] & = & 
     h_2[E,1,-E,-1,\dots] & = & \frac 12(E^2+1)\\[.05in]
    e_2[E,1,-E,-1,\dots] & = & 
     h_2[E,-1,-E,1,\dots] & = & \frac 12(E^2-1)\\[.05in]
   e_3[E,-1,-E,1,\dots] & = & 
     h_3[E,1,-E,-1,\dots] & = & \frac 16(E^3+E)\\[.05in]
   e_3[E,1,-E,-1,\dots] & = & 
     h_3[E,-1,-E,1,\dots] & = & \frac 16(E^3-5E)\\[.05in]
    e_4[E,-1,-E,1,\dots] & = & 
     h_4[E,1,-E,-1,\dots] & = & \frac{1}{24}(E^4-2E^2-3)\\[.05in]
     e_4[E,1,-E,-1,\dots] & = & 
     h_4[E,-1,-E,-1,\dots] & = & \frac{1}{24}(E^4-7E^2+9)\\[.05in]
   e_5[E,-1,-E,1,\dots] & = & 
     h_5[E,1,-E,-1,\dots] & = & \frac{1}{120}(E^5-10E^3-11)\\[.05in]
   e_5[E,1,-E,-1,\dots] & = & 
     h_5[E,-1,-E,1,\dots] & = & \frac{1}{120}(E^5-30E^4+89).
    \end{array} $$
\indent
It is easy to see (see equations~(\ref{eq:enhn1}), (\ref{eq:enhn2}),
(\ref{eq:enhn3}) below) that for all $i$ we have
  \beas e_i[E,0,-E,0,\dots] & = & h_i[E,0,-E,0,\dots]\\
      e_i[E,-1,-E,1,\dots] & = & h_i[E,1,-E,-1,\dots]\\
      e_i[E,1,-E,-1,\dots] & = & h_i[E,-1,-E,1,\dots]
   \eeas
These formulas, together with Theorem~\ref{thm:multi} and the
commutativity of the ring of symmetric functions, yield some results
about the equality of certain values of $N(\alpha,A,B)$. For instance,
if $n$ is odd, then $N(\alpha,A,B)$ depends only on the multiset of
parts of $\alpha$, not on their order, and also not on $A$ and $B$. If
$n$ is even, then $N(\alpha,A,B)$ depends only on the multiset of
parts of $\alpha$ and on which submultiset of these parts index the
elements of $A$ and $B$.

The specialization of $e_i$ and $h_i$ for small $i$ lead to some
nonumbral formulas for certain values of $N(\alpha,A,B)$. For
instance, let $k$ be odd, $\alpha=(3^k)$ (i.e., $k$ parts equal to 3),
$A=\emptyset$, so that $N((3^k),\emptyset,[k])$ is the number of
alternating permutations $a_1>a_2\leq a_3>a_4\leq a_5>\cdots\leq
a_{3k}$ (where $>$ and $\leq$ have their usual meaning) of the
multiset $\{1^3,2^3,\dots,k^3\}$. Then
  \beas N((3^k),\emptyset,[k]) & = & h_3^k[E,0,-E,0,\dots]\\
        & = & \frac{1}{6^k}E^k(E^2-2)^k\\ & = &
      \frac{1}{6^k}\sum_{j=0}^k{k\choose j}(-2)^{k-j}E_{2j+k}. \eeas
In the same way we obtain the formulas in \cite[pp.~265--266]{gessel}.

It is easy to find generating functions for the specializations of
$e_n$ and $h_n$ that we are considering, using the identities 
   \beas \sum_{n\geq 0} e_nt^n & = & \exp \sum_{j\geq 1} (-1)^{j-1}
     \frac{p_j}{j}\\
   \sum_{n\geq 0} h_nt^n & = & \exp \sum_{j\geq 1}
     \frac{p_j}{j}. \eeas
Namely,
  \bea \sum_{n\geq 0}e_n[E,0,-E,0,\dots]t^n  & = & \sum_{n\geq
     0}h_n[E,0,-E,0,\dots]t^n\nonumber\\ & = & 
   \exp E\tan^{-1}t \label{eq:enhn1}\\[.1in]
    \sum_{n\geq 0}e_n[E,1,-E,-1,\dots]t^n & = & 
    \sum_{n\geq 
     0}h_n[E,-1,-E,1,\dots]t^n \nonumber \\ & = &
      \displaystyle\frac{1}{\sqrt{1+t^2}}
     \exp E\tan^{-1}t \label{eq:enhn2}\\[.1in]
     \sum_{n\geq 0}e_n[E,-1,-E,1,\dots]t^n & = & \sum_{n\geq
     0}h_n[E,1,-E,-1,\dots] t^n\nonumber\\ & = & \sqrt{1+t^2}
     \exp E\tan^{-1}t \label{eq:enhn3} \eea
Equation~(\ref{eq:enhn1}) in fact is a restatement of (\ref{eq:eatt}).

\pagebreak

\end{document}